\documentclass[12pt,leqno,a4paper]{amsart}
\usepackage{amssymb,enumerate}
\usepackage{xcolor}
\usepackage{hyperref,backref}
\newcommand{\FF}{{\mathbb{F}}}
\newcommand{\QQ}{{\mathbb{Q}}}
\newcommand{\RR}{{\mathbb{R}}}

\newcommand{\bG}{{\mathbf{G}}}
\newcommand{\bH}{{\mathbf{H}}}
\newcommand{\bL}{{\mathbf{L}}}
\newcommand{\bS}{{\mathbf{S}}}
\newcommand{\bT}{{\mathbf{T}}}
\newcommand{\bZ}{{\mathbf{Z}}}

\newcommand{\subG}{\mathfrak{G}}
\newcommand{\fS}{{\mathfrak{S}}}

\newcommand{\cE}{{\mathcal{E}}}
\newcommand{\cN}{{\mathcal{N}}}

\newcommand{\cent}{\operatorname{C}}
\newcommand{\Dih}{\operatorname{Dih}}
\newcommand{\Ind}{{\operatorname{Ind}}}
\newcommand{\Irr}{{\operatorname{Irr}}}
\newcommand{\norm}{\operatorname{N}}
\newcommand{\St}{\operatorname{St}}
\newcommand{\Syl}{\operatorname{Syl}}
\newcommand{\type}{\operatorname}
\newcommand{\Uch}{{\operatorname{UCh}}}
\newcommand{\weyl}{\operatorname{W}}

\newcommand{\GL}{{\operatorname{GL}}}
\newcommand{\SL}{{\operatorname{SL}}}
\newcommand{\PSL}{{\operatorname{PSL}}}
\newcommand{\SU}{{\operatorname{SU}}}
\newcommand{\PGU}{{\operatorname{PGU}}}
\newcommand{\Sp}{{\operatorname{Sp}}}
\newcommand{\GO}{{\operatorname{GO}}}
\newcommand{\SO}{{\operatorname{SO}}}
\newcommand{\Spin}{{\operatorname{Spin}}}

\newcommand{\Chevie}{{\sf Chevie}{}}
\newcommand{\GAP}{{\sf{GAP}}}

\newcommand{\tw}[1]{{}^{#1}\!}
\newcommand{\wt}{\widetilde}
\newcommand{\RLG}{{R_\bL^\bG}}

\let\eps=\epsilon
\let\la=\lambda
\let\Om=\Omega
\let\tht=\theta
\let\vhi=\varphi
\let\ti=\times

\theoremstyle{theorem}
\newtheorem{thm}{Theorem}[section]
\newtheorem{lem}[thm]{Lemma}
\newtheorem{prop}[thm]{Proposition}

\newtheorem{thmA}{Theorem}
\newtheorem{corA}[thmA]{Corollary}
\newtheorem{conjA}[thmA]{Conjecture}
\newtheorem{definition}[thmA]{Definition}

\theoremstyle{definition}
\newtheorem{rem}[thm]{Remark}

\makeatletter
\def\namedlabel#1#2{\begingroup #2%
    \def\@currentlabel{#2}%
    \phantomsection\label{#1}\endgroup
}

\begin{document}

\title[Picky Conjecture]{The Picky Conjecture for groups of Lie type\\
  in non-defining characteristic}

\author{Gunter Malle}
\address[G. Malle]{FB Mathematik, RPTU Kaiserslautern, Postfach 3049,
  67653 Kaisers\-lautern, Germany.}
\email{malle@mathematik.uni-kl.de}

\author{A. A. Schaeffer Fry}
\address[A. A. Schaeffer Fry]{Dept. Mathematics - University of Denver, Denver,
  CO 80210, USA}
\email{mandi.schaefferfry@du.edu}

\begin{abstract}
Recently, Moret\'o announced a conjecture, known as the Picky Conjecture,
proposing new character correspondences for finite groups that extend the
McKay Conjecture in certain situations. We prove the Picky Conjecture for all
quasi-simple groups of Lie
type for non-defining primes. In favourable situations, we also obtain the
stronger version postulating preservation of character values up to sign, and
we show this stronger version holds in general when assuming certain natural
properties of Lusztig's Jordan decomposition. Along the way, we complete the
determination of picky elements in these groups by classifying semisimple
picky $p$-elements for $p\in\{2,3\}$.
\end{abstract}



\thanks{The authors thank the anonymous referees for their careful reading and
 constructive comments. The first author gratefully acknowledges financial
 support by the DFG -- Project-ID 286237555. The second author is grateful for
 support of a CAREER grant through the U.S. National Science Foundation,
 Award No. DMS-2439897.}

\date{\today}

\maketitle


\section{Introduction}   \label{sec:intro}

The investigation of global to local character correspondences for finite
groups has been a very fruitful area of research in the past several decades.
Here we investigate a new type of correspondence that has recently been
proposed by Alex Moret{\'o} at a conference in Val\`encia in January 2025,
coined the \emph{Picky Conjecture}, and prove it for quasi-simple groups of
Lie type in non-defining characteristic.

Let $G$ be a finite group and $\ell$ a prime. We say that an $\ell$-element
$x\in G$ is \emph{picky} if it lies in a unique Sylow $\ell$-subgroup $P$
of~$G$. For any $x\in G$, we let $\Irr^x(G)$ denote the set of complex
irreducible characters of $G$ that take a non-zero value at $x$.

We are concerned with the following recent conjecture \cite[Conj.~A]{M26}
(see also \cite{MR25}, where Moret{\'o} and Noelia Rizo began studying the
conjecture in some detail).

\begin{conjA}[{Picky Conjecture}]   \label{conj:AN}
 Let $G$ be a finite group, $\ell$ a prime, $P\in\Syl_\ell(G)$ and $x\in P$ a
 picky $\ell$-element. Then there exists a bijection
 $$f:\Irr^x(G)\to\Irr^x(\norm_G(P))$$
 such that for each $\chi\in\Irr^x(G)$,
 \begin{itemize}
  \item[\namedlabel{deg}{(1)}] $\chi(1)_\ell=f(\chi)(1)_\ell$, and
  \item[\namedlabel{fields}{(2)}] $\QQ(\chi(x))=\QQ(f(\chi)(x))$.
 \end{itemize}
  We will say the \emph{Picky Conjecture holds for} $(G,\ell,x)$ if \ref{deg}
  and \ref{fields} hold.
  \end{conjA}

As observed in \cite{M26}, the Picky Conjecture has strong connections
with the recently-proved McKay Conjecture, as well as with its Galois
equivariant version proposed by Navarro. In particular, if $G$ has some picky
$\ell$-element, then the Picky Conjecture at the prime~$\ell$ implies the McKay
Conjecture for~$G$ at~$\ell$. 

We will work with two stronger versions of the Picky Conjecture.

\begin{definition}[Picky+ Conjecture]   \label{def:+}
 We will say that the \emph{Picky+ Conjecture holds for} $(G, \ell, x)$ if the
 Picky Conjecture holds for $(G, \ell, x)$ and the following further holds:
 \begin{itemize}
  \item[\namedlabel{ellpart}{(3)}] $\chi(x)_\ell=f(\chi)(x)_\ell\quad$ for all $\chi\in\Irr^x(G)$.
 \end{itemize}
\end{definition}

Here $\chi(x)_\ell$ denotes the $\ell$-adic valuation in $\RR_{\geq 0}$:
namely, for $\alpha\in\overline{\QQ}$, we write
$\alpha_\ell:=\left|\cN(\alpha)\right|_\ell^{1/[\QQ(\alpha):\QQ]}$ where
$\cN:=\cN_{\QQ(\alpha)/\QQ}$ is the norm map.

\begin{definition}[Strong Picky Conjecture]   \label{def:strong}
 We will say that the \emph{Strong Picky Conjecture holds for} $(G,\ell,x)$ if
 the Picky Conjecture holds and the following strong version of~\ref{fields}
 holds (which implies \ref{fields} and~\ref{ellpart}):
 \begin{itemize}
  \item[\namedlabel{values}{(2*)}] $\chi(x)=\pm f(\chi)(x)\quad$ for all $\chi\in\Irr^x(G)$.
 \end{itemize}
\end{definition}

We remark that the Strong Picky Conjecture does not hold in general. Indeed,
as pointed out in \cite{M26}, the groups in \cite[Prop.~1.3]{blaumichler} with non-cyclic TI Sylow
$\ell$-subgroups provide counterexamples to~(2*)
 (see \cite[Thm.~3.1]{M26}). 
However, Moret{\'o}--Navarro--Rizo recently showed the Strong Picky Conjecture
for $\ell$-solvable groups for odd primes $\ell$ in \cite[Thm~A]{MNR}, and
Mart{\'i}nez Madrid recently showed in \cite[Thm~A]{MartinezMadrid} that it
holds for symmetric groups. In this paper, we find that the Strong Picky
Conjecture holds for groups of Lie type
in non-defining characteristic, assuming some standard conjectures on
characters of groups of Lie type (that have been shown to hold in many cases;
see Proposition~\ref{prop:strongpickyreg} below).

Let $G=\bG^F$ be a finite group of Lie type, where $\bG$ is a connected
reductive linear algebraic group in characteristic~$p>0$ and $F:\bG\to\bG$ is a
Steinberg map.
In \cite{Ma25}, the first author began the classification of picky elements for
such $G$, as well as the larger study of subnormalisers of $\ell$-elements and
the corresponding more general version of Conjecture \ref{conj:AN} for
arbitrary $\ell$-elements. It turns out that for $\ell=p$, there exist picky
$\ell$-elements even in highly non-abelian Sylow subgroups and in arbitrary
rank; the character tables of Sylow normalisers in this case are not
well-understood.

In contrast, if $p\ne\ell>3$ and $\bG$ is simple of simply connected type,
by \cite[Thm~5.9]{Ma25} then $G$ has no picky $\ell$-elements unless $P$ is
abelian. Hence, our first focus in this paper is on the situation that
$\ell\neq p$ and that $P$ is abelian. In particular, we show that in this
situation, the Picky Conjecture holds (see Proposition~\ref{prop:nonreg}).
Combining this with an investigation of the finitely many groups with picky
3-elements we obtain:

\begin{thmA}   \label{thm:pickyodd}
 Let $G=\bG^F$ be a finite group of Lie type, where $\bG$ is simple and simply
 connected in characteristic $p$ and $F$ is a Steinberg morphism. Then the
 Picky+ Conjecture holds for $G$ for any prime $\ell\geq 3$ different from~$p$.
\end{thmA}

In our proof of Theorem \ref{thm:pickyodd}, in addition to the work in
\cite{Ma25}, we also rely on the McKay bijection proved by the first author and
Sp{\"a}th in \cite{Ma07, spath10}. 

We also classify the picky semisimple $\ell$-elements in the case that
$\ell=2$, and we obtain the stronger version in this case:

\begin{thmA}   \label{thm:pickyeven}
 Let $G=\bG^F$ be a finite group of Lie type, where $\bG$ is simple and simply
 connected in characteristic $p\neq 2$ and $F$ is a Steinberg morphism. Then the
 Strong Picky Conjecture holds for $G$ for the prime $\ell=2$.
\end{thmA}

With this we also eventually obtain the Picky Conjecture for quasi-simple
groups of Lie type in non-defining characteristic:

\begin{corA}   \label{cor:quasisimple}
 The Picky Conjecture holds for any quasi-simple group~$H$ such that
 $H/\bZ(H)$ is simple of Lie type in characteristic~$p$ for any
 prime~$\ell\ne p$.
\end{corA}

The paper is organised as follows. In Section~\ref{sec:abelian}, we treat the
generic case of finite reductive groups with abelian Sylow $\ell$-subgroups,
and also discuss the Suzuki and Ree groups. In Section~\ref{sec:picky 2}, we
determine quasi-simple groups of Lie type with semisimple picky 2-elements,
and in the final section, Section~\ref{sec:23exceptcover}, we complete the
proofs of our main Theorems~\ref{thm:pickyodd} and~\ref{thm:pickyeven} and of
Corollary~\ref{cor:quasisimple} by first classifying and then dealing with
picky 2- and 3-elements, as well as groups with exceptional covering groups.
\medskip

We remark that shortly after acceptance of this article, Damiano Rossi
announced a counter-example to the Picky Conjecture in the case of non-abelian
Sylow subgroups.

\section{The case of abelian Sylow subgroups}   \label{sec:abelian}
In this section, we prove Theorem \ref{thm:pickyodd} when the Sylow subgroups
of $G$ are abelian. First, let $G$ be an arbitrary finite group with abelian
Sylow subgroups. In this situation, the following reduces the study of the Picky
Conjecture to a strengthening of bijections in the McKay Conjecture. For this
recall that $\Irr_{\ell'}(G)$ denotes the set of complex irreducible characters
of $G$ of degree prime to~$\ell$.

\begin{lem}[Prop.~4.5 of \cite{M26}]   \label{lem:abelian}
 Let $\ell$ be a prime and suppose that $G$ is a finite group with abelian
 Sylow $\ell$-subgroups. Let $x\in G$ be a picky $\ell$-element. Then
 $\Irr^x(G)=\Irr_{\ell'}(G)$.
\end{lem}

\begin{rem}\label{rem:picky+}
Note that in the situation of Lemma \ref{lem:abelian}, \ref{ellpart} of
Definition~\ref{def:+} follows by general properties of character values.
Indeed, let $\chi\in\Irr_{\ell'}(G)$. Now $\chi(x)$ is a sum of $n:=\chi(1)$
many $\ell$-power order roots of unity, and if $m$ denotes the length of its
Galois orbit, then the norm $\cN(\chi(x))$ is an
integer that is a sum of $n^m$ $\ell$-power roots of unity. Now for $a\ge2$
the sum over all $\ell^{a-1}(\ell-1)$ primitive $\ell^a$th roots of unity is
zero (inductively, since the sum over all $\ell^a$th roots of unity vanishes),
so we see that $\cN(\chi(x))$ equals a sum of $k$ $\ell$th roots of unity where
$k\equiv n^m\pmod\ell$. Finally, the sum over the $\ell-1$ primitive $\ell$th
roots of unity equals $-1$, so $\cN(\chi(x))\equiv n^m\pmod\ell$ is prime to
$\ell$. This of course also holds for $\norm_G(P)$ in place of $G$, where
$P\in\Syl_\ell(G)$.
\end{rem}

As the McKay Conjecture has been proven (see \cite{CS24}), we see from
Lemma~\ref{lem:abelian} and Remark~\ref{rem:picky+} that to prove the Picky+
Conjecture, resp.\ the Strong Picky Conjecture, for $(G,\ell, x)$ when
$P\in\Syl_\ell(G)$ is abelian, it suffices to prove item \ref{fields} of
Conjecture~\ref{conj:AN}, resp.\ \ref{values} of Definition~\ref{def:strong},
for some known McKay bijection
$f\colon \Irr_{\ell'}(G)\to \Irr_{\ell'}(\norm_G(P))$.

We also remark that the (Strong) Picky Conjecture is known for the case that
$P$ is cyclic:

\begin{lem}[Rizo]   \label{lem:cyclic}
 Let $\ell$ be a prime and suppose that $G$ has cyclic Sylow $\ell$-subgroups.
 Then the Strong Picky Conjecture holds for all picky $\ell$-elements of~$G$.
\end{lem}

We note that in \cite[Thm~4.7]{M26}, Lemma~\ref{lem:cyclic} is proven in the
case that Sylow $\ell$-subgroups have order $\ell$, and it is noted in the
subsequent comments that the arbitrary cyclic case was proved by 
N.~Rizo, using Dade's work on cyclic blocks. 

Before setting the stage for our main results, we observe the following,
reducing the conjectures at hand for quasi-simple groups to the case of
universal covering groups. For this let us say that a character bijection $f$
for $G$ as in our conjectures \emph{preserves central characters} if $f(\chi)$
lies over the same character of $\bZ(G)$ as $\chi$ for each $\chi\in\Irr^x(G)$.

\begin{lem}   \label{lem:universalcoversuffices}
 Let $S$ be a finite non-abelian simple group with universal covering group
 $G$. Suppose that the Picky (resp.~Picky +, resp. Strong Picky) conjecture
 holds for $(G,\ell,x)$ for all picky $\ell$-elements $x$ of $G$ with
 a bijection preserving central characters.
 Then given a quasi-simple group $H$ with $H/\bZ(H)=S$, the Picky (resp.~Picky
 +, resp. Strong Picky) conjecture holds for $(H,\ell,y)$ for all picky
 $\ell$-elements $y$ in $H$.
\end{lem}

\begin{proof}
We may write $H=G/Z$ for some central subgroup $Z\leq \bZ(G)$. Let $y\in H$ be
a picky $\ell$-element lying in the Sylow $\ell$-subgroup $\bar{P}$ of $H$.
By \cite[Lem.~2.4]{Ma25}, $y=xZ$ for some picky $\ell$-element $x$ in $G$.
Write $\bar{P}=PZ/Z$ for some $P\in\Syl_\ell(G)$, and note that necessarily
$x\in P$. We may view any character of $H$ as a character of $G$ by inflation,
and as such, we have $\Irr^y(H)=\{\chi\in\Irr^x(G)\mid Z\leq \ker\chi\}$, and
the analogous statement holds for
$\Irr^y(\norm_H(\bar{P}))\subseteq \Irr^x(\norm_G(P))$. From here, the
additional assumption that $f$ preserves central characters yields the result.
\end{proof}

\subsection{The Setting}   \label{subsec:setting}
We consider the following setting, where we refer the reader to \cite{MT11} and
\cite{GM20} for notation and background on linear algebraic groups and the
character theory of finite groups of Lie type. As before, let $\bG$ be simple
of simply connected type in characteristic~$p$ with a Frobenius map
$F:\bG\to\bG$ with respect to an $\FF_q$-structure and $G:=\bG^F$.
Let $(\bG^\ast, F)$ be dual to $(\bG, F)$, and write $G^\ast:=(\bG^\ast)^{F}$. 

We let $\ell\ne p$ be a prime and suppose that $P\in\Syl_\ell(G)$ is abelian.
This means that $\ell$ satisfies the condition ($\ast$) of
\cite[Sec.~2]{Malle14}. Further, this means that $\ell$ is an odd prime good
for $\bG$ which divides neither $[\bZ(\bG)^F:\bZ^\circ(\bG)^F]$ nor
$[\bZ(\bG^\ast)^F:\bZ^\circ(\bG^\ast)^F]$, and that further $\ell\geq 5$ if
$G$ equals $\tw{3}\type{D}_4(q)$. (See \cite[Lemma~2.1, Prop.~2.2]{Malle14}.)

Let $d:=d_\ell(q)$ be the order of $q$ modulo $\ell$ and let $\bS$ be a Sylow
$d$-torus of $(\bG, F)$ containing $P$, which exists by our assumptions on
$\ell$ (see the proof of \cite[Prop.~2.3]{Malle14}). Let $\bL:=\cent_\bG(\bS)$
be the corresponding minimal $d$-split Levi subgroup of $\bG$ and write
$L:=\bL^F$. We further let $\bL^\ast:=\cent_{\bG^\ast}(\bS^\ast)$, where
$\bS^\ast$ is a Sylow $d$-torus for $(\bG^\ast, F)$ dual to $\bS$ and (thus)
$\bL^\ast$ is dual to~$\bL$. (Such a choice exists by \cite[Lemma~3.3]{Ma07}.)
We then let $L^\ast:={\bL^\ast}^F$ and also write
$N^\ast:=\norm_{G^\ast}(\bS^\ast)$.

From the first author's work in \cite{Malle14}, we see:

\begin{lem}   \label{lem:strab}
 In the above situation, we have
 \begin{enumerate}[\rm(a)]
  \item $\norm_G(P)=\norm_\bG(\bL)^F$;
  \item $L=P\times \operatorname{O}^\ell(L)$;
  \item $\ell\nmid|\weyl_G(\bL)|$, where $\weyl_G(\bL):=\norm_\bG(\bL)^F/L$; 
  \item any picky $\ell$-element $x\in G$ lies in a unique $G$-conjugate of $L$;
   and
  \item if $x\in L$ is a picky $\ell$-element, then $x^g\in L$ if and only if
   $g\in\norm_G(P)$.
\end{enumerate}
\end{lem}

\begin{proof}
The first three items are from \cite[Prop.~2.4]{Malle14}.  The latter two items
follow, since if $x\in P$ is picky and $x\in L^g$ for some $g\in G$, then
$x\in P^g$ by (b). Hence $g\in\norm_G(P)=\norm_\bG(\bL)^F$ by (a). Then
$L^g= (\bL\cap G)^g=\bL\cap G=L$.
\end{proof}

Further, the picky elements as above were classified in
\cite[Thm~5.2]{Ma25}:

\begin{thm}   \label{thm:gunter}
 In the above situation, $x\in P$ is picky if and only if $\cent_\bG(x)=\bL$.
\end{thm}

Write $N:=\norm_G(P)=\norm_\bG(\bL)^F$. In what follows, we will consider a
bijection 
\[\Om\colon \Irr_{\ell'}(G)\to\Irr_{\ell'}(N)\]
as constructed in \cite[Thm~7.8]{Ma07}, which uses the results of \cite{spath09,
spath10} to obtain the local extensions that yield this bijection. Throughout,
$\Om$ will continue to denote such a bijection, and we will argue that any
such $\Om$ satisfies the properties for $f$ postulated in the Picky Conjecture. 

We also note here that by \cite[Thm~7.8(c)]{Ma07}, $\Om$ has the
additional property that $\chi$ and $\Om(\chi)$ lie above the same character
of $\bZ(G)$ for any $\chi\in\Irr_{\ell'}(G)$. Now, if $S$ is a simple group
of Lie type in characteristic $p\neq\ell$ with non-exceptional Schur multiplier
but not of Suzuki or Ree type, then we know that its universal covering group
is of the form $G=\bG^F$ above.
Hence, together with the results of \cite{Ma25}, in such a case, to prove the
Picky+ Conjecture for all quasi-simple groups $H$ with $H/\bZ(H)\cong S$ and
primes $\ell>3$ with $\ell\neq p$, it suffices to prove that
$\QQ(\chi(x))=\QQ(\Om(\chi)(x))$ where $G$, $x$, $\Om$ are as above. (Note that
the cases $\ell\in\{2, 3\}$ or $S$ has an exceptional Schur multiplier are
treated in Section~\ref{sec:23exceptcover} below, and the Suzuki and Ree groups
are discussed in Section~\ref{subsec:suzree} below.)

\subsection{Values on the local side}   \label{subsec:local}
Keep the setting from Section \ref{subsec:setting}. We first establish the
relevant values for characters in $\Irr_{\ell'}(N)$. Let $x\in P$ be a picky
$\ell$-element. By Lemma~\ref{lem:abelian} we have $\Irr^x(G)=\Irr_{\ell'}(G)$.

So let $\chi\in\Irr_{\ell'}(G)$. 
In this situation, the character $\Om(\chi)\in\Irr_{\ell'}(N)$ is of the form 
\[\Om(\chi)=\Ind_{N_\tht}^N(\hat\tht \eta),\]
where $\hat\tht\in\Irr(N_\tht)$ is an extension of some character
$\tht\in\Irr_{\ell'}(L)$ (whose existence is guaranteed by
\cite[Thm~1.1]{spath10}) and $\eta$ is a character of
$N_\tht/L=\weyl_G(\bL)_\tht$ viewed as a character of $N_\tht$ by inflation.
We will say that $\Om(\chi)$ corresponds to the pair $(\theta, \eta)$. Note
here that since $P\lhd L$ is abelian, in fact $\Irr_{\ell'}(L)=\Irr(L)$ by
Clifford theory.

Then since $x^g\in L\leq N_\tht$ for any $g\in N$ and $\eta$ is trivial on
$L$, we have
\[\Om(\chi)(x)=\frac{1}{|N_\tht|}\sum_{g\in N; x^g\in N_\tht} \hat{\tht}(x^g)\eta(x^g)
  =\frac{1}{|N_\tht|}\sum_{g\in N} {\tht}(x^g)\eta(x^g)
  =\frac{\eta(1)}{|N_\tht|}\sum_{g\in N; x^g\in L} {\tht}(x^g).\]
Thus we see 
\begin{equation}\label{eq:local}
\Om(\chi)(x)=\eta(1)\frac{|L|}{|N_\tht|}\Ind_{L}^N(\tht) (x)
\end{equation}
and hence in particular
\begin{equation}\label{eq:local2}
\QQ(\Om(\chi)(x))=\QQ(\Ind_{L}^N(\tht) (x)).
\end{equation}

\subsection{Values on the global side}   \label{subsec:global}
We continue to keep the setting from the previous sections and now investigate
the relevant values for $\chi\in\Irr_{\ell'}(G)$.  Recall that
$\cent^\circ_\bG(x)=\cent_\bG(x)=\bL$ by Theorem~\ref{thm:gunter}. By
the Curtis--Schewe formula (\cite[Cor.~3.3.13]{GM20} and the discussion after),
we have:
\begin{equation}\label{eq:valchi}
  \chi(x)=\tw{\ast}\RLG(\chi)(x)
         =\sum_{\vhi\in\Irr(L)} \langle \chi, \RLG(\vhi)\rangle \vhi(x),
\end{equation}
where $\RLG$ and $\tw{\ast}\RLG$ denote Lusztig induction, resp.\ 
restriction.

Now, if $\chi$ lies in the rational Lusztig series $\cE(G,s)$ of $G$
corresponding to a semisimple element $s\in G^\ast$, then up to
$G^\ast$-conjugacy, we may assume $s\in L^\ast$ as $\chi\in\Irr_{\ell'}(G)$.
(See \cite[Prop.~7.2]{Ma07}.) From the work constructing the bijection~$\Om$
in \cite[Thm~7.5, 7.8 and Prop.~7.7]{Ma07}, 
we see we may assume further that $\tht\in\cE(L, s)$ if $\chi$ is such
that $\Om(\chi)$ corresponds to $(\tht,\eta)$ as in Section~\ref{subsec:local}. 

\begin{rem}   \label{rem:classicalcent}
In the above situation, we note that since $\bL^*$ is minimal $d$-split,
$\cent_{\bL^\ast}(s)$ contains only
components of classical type $\type{A}$, $\type{B}$, $\type{C}$, $\type{D}$.
This is clear if $d=d_\ell(q)$ is a so-called \emph{regular number} for
$(\bG, F)$, since then $\bL^\ast$ is a torus (see \cite[Ex.~3.5.7]{GM20}).
If $d$ is not a regular number, the claim follows from \cite[Tab.~3.3]{GM20}
and \cite[Tab.~1]{BMM93} by inspection.
\end{rem}

In the following, we consider a regular embedding
$\iota\colon\bG\hookrightarrow \wt\bG=\bG\bZ(\wt\bG)$ and its dual
epimorphism $\iota^\ast\colon \wt\bG^\ast\to\bG^\ast$, as in
\cite[Prop.~1.7.5]{GM20}, and extend $F$ naturally to $\wt\bG$, resp.\ to
$\wt\bG^\ast$. We also recall that, thanks to Lemma \ref{lem:cyclic}, we are
particularly interested in the case that $P$ is not cyclic.

\begin{prop}   \label{prop:dHCseries}
 Keep the setting of Section \ref{subsec:setting}. Thus, $G=\bG^F$ where
 $\bG$ is simple of simply connected type in characteristic $p\neq \ell$, with
 $F\colon \bG\rightarrow\bG$ a Frobenius map; $P\in\Syl_\ell(G)$ is abelian
 lying in a Sylow $d$-torus $\bS$ with $d:=d_\ell(q)$; and $L:=\bL^F$, where
 $\bL:=\cent_\bG(\bS)$. Further, write
 $\wt{\bL}:=\bL\bZ(\wt\bG)$ and $\wt{L}:=\wt\bL^F$.

 Let $\chi\in\Irr_{\ell'}(G)\cap\cE(G,s)$ such that $\Om(\chi)$ corresponds to
 $(\tht, \eta)$ as in Section \ref{subsec:local}, where $s\in L^\ast$ and
 $\tht\in\cE(L,s)$. If
 \begin{enumerate}
  \item[\rm(1)] $J_s^\bG\circ \RLG
     =R_{\cent_{\bL^\ast}(s)}^{\cent_{\bG^\ast}(s)}\circ J_s^\bL$; or
  \item[\rm(2)] $P$ is non-cyclic,
 \end{enumerate}
 then $\langle\chi,\RLG(\tht')\rangle\neq 0$ for some $\wt{L}$-conjugate
 $\tht'$ of $\tht$. In Case~(1), if $\cent_{\bL^\ast}(s)$ is further connected,
 then we may take $\tht'=\tht$.
\end{prop}

\begin{rem}   \label{rem:regdHC}
 We note that if $d$ is a regular number for $G$, then (1) holds and
 $\cent_{\bL^\ast}(s)$ is connected, so that $\theta'=\theta$ in the situation
 of Proposition \ref{prop:dHCseries}. Indeed, then
 $\cent_{\bL^\ast}(s)=\bL^\ast$ is a torus, hence connected, and further the
 commutation holds by the defining properties of Jordan decomposition
\cite[Thm~2.6.22]{GM20}.
\end{rem}

\begin{proof}[Proof of Proposition \ref{prop:dHCseries}]
From \cite[Prop.~7.7]{Ma07}, $\tht$ corresponds under some Jordan decomposition
$J_s^\bL$ to a character $J_s^\bL(\theta)$ of $\cent_{L^\ast}(s)$
lying above some unipotent character
$\la\in\Uch(\cent^\circ_{\bL^\ast}(s)^F):=\cE(\cent^\circ_{\bL^\ast}(s)^F,1)$
and furthermore 
\begin{equation}\label{eq:relweyl}
\weyl_{\cent_{G^\ast}(s)}(\cent^\circ_{\bL^\ast}(s))_\la=\weyl_{N^\ast}(s,\la)
  \cong \weyl_G(\bL)_\tht
  \end{equation}
under duality, where we recall that $N^\ast:=\norm_{G^\ast}(\bS^\ast)$ as in
Section~\ref{subsec:setting}. Also, here
$\cent_{L^\ast}(s)=\cent_{\cent_{\bG^\ast}(s)}(\bS)^F$ and
$\cent_{L^\ast}^\circ(s):=\cent_{\bL^\ast}^\circ(s)^F
  =\cent_{\cent^\circ_{\bG^\ast}(s)}(\bS)^F$, and the latter is
also a minimal $d$-split Levi subgroup of $(\cent_{\bG^\ast}^\circ(s), F)$
according to the proof of \cite[Cor.~6.6]{Ma07}. 

On the local side, we have seen that $\Om(\chi)$ corresponds to the pair
$(\tht, \eta)$, or the triple $(s,\la,\eta)$.
On the global side, $\chi$ also corresponds to the triple $(s,\la,\eta)$, which
in the context of \cite[Sec.~6,7]{Ma07} means that $\chi$ corresponds via some
Jordan decomposition $J_s^\bG$ to
$J_s^\bG(\chi)\in\Uch(\cent_{G^\ast}(s))$ lying in the $d$-Harish-Chandra
series of $(\cent_{L^\ast}^\circ(s), \la)$ as defined in \cite[Thm~4.6]{Ma07}.
(Here $J_s^\bG(\chi)$, in turn, corresponds to $\eta$ under the bijection of
\cite[Thm~4.6(b)]{Ma07}.) Then in particular, $J_s^\bG(\chi)$ is a constituent
of
\begin{equation}\label{eq:IndRLG}
\Ind_{\cent_{G^\ast}^\circ(s)}^{\cent_{G^\ast}(s)}(R_{\cent_{\bL^\ast}^\circ(s)}^{\cent_{\bG^\ast}^\circ(s)}(\la))
  =R_{\cent_{\bL^\ast}(s)}^{\cent_{\bG^\ast}(s)}(\Ind_{\cent_{L^\ast}^\circ(s)}^{\cent_{L^\ast}(s)}(\la)),
\end{equation}
where $R_{\cent_{\bL^\ast}(s)}^{\cent_{\bG^\ast}(s)}$ is as defined in
\cite[Def.~4.8.8]{GM20} and the equality follows from \cite[Cor.~2.4(ii)]{DM94}
by adjunction. Note further that by \cite[3.3.20]{GM20}, the constituents of
$\RLG(\tht)$ lie in $\cE(G,s)$.

Assume first that (1) holds. If further $\cent_{\bL^\ast}(s)$ is connected,
then $\chi$ lies in the $d$-Harish-Chandra series of $(\bL,\tht)$. Indeed,
we then have $J_s^\bL(\tht)=\la$ and
$J_s^\bG(\RLG(\tht))=R_{\cent_{\bL^\ast}(s)}^{\cent_{\bG^\ast}(s)}(\la)$
contains $J_s^\bG(\chi)$. Since $J_s^\bG$ is a bijection from $\cE(G, s)$ to
$\Uch(\cent_{G^\ast}(s))$, this means that $\langle\chi,\RLG(\tht)\rangle\neq0$.
If $\cent_{\bL^\ast}(s)$ is disconnected, $J_s^\bG(\chi)$ is a
constituent of $R_{\cent_{\bL^\ast}(s)}^{\cent_{\bG^\ast}(s)}(\hat\la)$, for
some character~$\hat\la$ of $\cent_{G^*}(s)$ lying above $\la$.
Let $\tht'\in\cE(L, s)$ be such that
$J_s^\bL(\tht')=\hat\la$. Then since $J_s^\bL(\tht')$ and $J_s^\bL(\tht)$ lie
above the same unipotent character of $\cent_{L^\ast}^\circ(s)$, it follows
from the analysis in \cite[Sec.~2.6]{GM20} that $\tht$ and $\tht'$ are
$\wt{L}$-conjugate. Here $J_s^\bG(\chi)$ is a constituent of
$J_s^\bG(\RLG(\tht'))$, and again we see $\langle\chi,\RLG(\tht')\rangle\neq 0$.

In Remark \ref{rem:regdHC} we see that if $d$ is regular, we have (1) and
$\cent_{\bL^\ast}(s)$ is connected, so $\chi$ lies in the $d$-Harish-Chandra
series of $(\bL, \theta)$. In particular, we may assume now that $d$ is not a
regular number.

Now assume that (2) holds. If $G$ is of exceptional type, since $d$ is not
regular and $P$ is not cyclic by~(2), we conclude that $G=\type{E}_7(q)$ and
$d=4$. Hence, we may assume that either $G$ is of classical type or that
$G=\type{E}_7(q)$ with $d=4$, and recall from Remark \ref{rem:classicalcent}
that $\cent_{\bL^\ast}(s)$ contains only components of classical type in either
situation. Now, since $\bZ(\wt\bG)$ is connected, the Mackey formula holds
for $\wt\bG$ by \cite[Thm~3.3.7(5)]{GM20}. Then further
$J_{\wt s}^{\wt\bG}\circ R_{\wt\bL}^{\wt\bG}
     =R_{\cent_{\wt\bL^\ast}(\wt s)}^{\cent_{\wt\bG^\ast}(\wt s)}\circ J_{\wt s}^{\wt \bL}$ 
for $\wt{s}\in\wt\bG^{\ast F}$ with $\iota^\ast(\wt{s})=s$, using
\cite[Thm~4.7.2]{GM20} when $G$ is of classical type and using the reasoning
in the second two paragraphs of the proof of \cite[Thm~4.7.5]{GM20} if
$G=\type{E}_7(q)$. 

Note here that $\wt\bL$ is a minimal $d$-split Levi subgroup of $\wt\bG$.
Further, since $\ell\nmid |\wt{G}/G\bZ(\wt{G})|$ by our assumptions, we see
that $\Irr_{\ell'}(\wt{G})$ is exactly the set of characters lying above those
in $\Irr_{\ell'}(G)$. Now, in the current situations, each character of
$\Irr(\wt{L})$ extends to its inertia group in $\wt{N}:=\norm_{\wt{G}}(\bS)$.
Indeed, for types $\type{A}$, $\type{B}$, $\type{C}$, and $\type{D}$, this is
from \cite[Cor.~5.14]{CS17a}, \cite[Sec.~5.5]{CS19}, \cite[Thm~6.1]{CS17b},
and \cite[Thms~4.1, ~6.9]{CS24}, respectively. (See also
\cite[Prop.~3.20]{MS16} for the case $d\in\{1,2\}$.) For type $\type{E}_7$,
this is shown in the course of checking the conditions for \cite[Thm~4.2]{CS19}
in the proof of \cite[Thm~A]{CS19}. In particular, there is an extension map
$\wt\Lambda$ with respect to $\wt{L}\lhd \wt{N}$, and hence we may apply the
parametrisation of $\Irr_{\ell'}(\wt{G})$ from \cite[Sec.~4]{CS13}. With this,
considerations similar to above, now using this parametrisation, yield that we
may find a character
$\wt\chi\in\cE(\wt{G}, \wt{s})\cap \Irr_{\ell'}(\wt{G})$ lying above $\chi$ 
such that $\langle\wt\chi, R_{\wt\bL}^{\wt\bG} (\wt{\tht})\rangle\neq 0$, where
$\wt\tht\in\cE(\wt{L}, \wt{s})$ lies above $\tht$. Then by
\cite[Cor.~3.3.25]{GM20}, we have $\langle\chi, \RLG(\tht')\rangle\neq 0$ for
some $\wt{L}$-conjugate $\tht'$ of~$\tht$.
\end{proof}

The following proves Theorem \ref{thm:pickyodd} for $\ell>3$, thanks to
\cite[Thm~5.9]{Ma25}.

\begin{prop}   \label{prop:nonreg}
 Keep the setting in Section \ref{subsec:setting}. Then
 Conjecture~\ref{conj:AN}\ref{fields} holds. That is, the Picky+ Conjecture
 holds for $(G,\ell,x)$ for all picky $\ell$-elements $x\in G=\bG^F$ when
 $\ell$ is an odd prime such that the Sylow $\ell$-subgroups of $G$ are abelian.
\end{prop}

\begin{proof}
Recall from Remark \ref{rem:picky+} that the second statement follows from the
first. By Lemma \ref{lem:cyclic} we may assume that $P\in\Syl_\ell(G)$ is
non-cyclic.
Hence we are as in the hypotheses in Proposition~\ref{prop:dHCseries}(2),
so $\langle \chi, \RLG(\tht')\rangle\neq 0$ for some $\wt{L}$-conjugate
$\tht'$ of $\tht$, where we recall that $\chi\in\Irr_{\ell'}(G)\cap\cE(G,s)$
with $s\in L^\ast$ and $\tht\in\cE(L,s)$.

Also, recall that, without loss, $x\in P\leq \bS$. Note that since $x\in \bS$
is centralised by any $g\in\wt{L}$, we see \eqref{eq:local2} holds with $\tht$
replaced with $\tht'$. Then for the purposes of proving that
$\QQ(\chi(x))=\QQ(\Om(\chi)(x))$, we may assume that
$\langle\chi,\RLG(\tht)\rangle\neq 0$.

To prove that $\QQ(\chi(x))=\QQ(\Om(\chi)(x))$, \eqref{eq:local2} and
\eqref{eq:valchi}, together with Lemma~\ref{lem:strab}(e), yield that it now
suffices to show that
$$\sum_{\vhi\in\Irr(L)} \langle \chi, \RLG(\vhi)\rangle \vhi(x) 
 =r\sum_{g\in N/N_\tht} \tht^g(x)\quad\text{for some $r\in\QQ^\times$}.$$

Thanks to \cite[Lemma~3.3, Prop.~7.7]{Ma07}, we have that
$N/N_\tht\cong N^\ast/L^\ast\cent_{N^\ast}(s)_\la$ where $\chi$ corresponds
to the triple $(s,\la,\eta)$ as before. Using \cite[Prop.~1.9]{CE99}, we see
that for $g\in N/N_\tht$, we have $\cE(L, s)^g=\cE(L, s^{g^\ast})$ where
$g\mapsto g^\ast$ denotes some such isomorphism.
Further, by \cite[Lemma~2.2]{Ma07}, each character in $\cE(L, s)^g$ lies above
the same character of $\bZ(L)$, so that $\vhi(x)/\vhi(1)=\tht^g(x)/\tht(1)$
is constant for $\vhi\in\cE(L, s)^g$ since $\tht\in\cE(L, s)$ and
$x\in P\leq\bZ(L)$.
 
Now, let $\vhi\in \Irr(L)$ with $\langle\chi, \RLG(\vhi)\rangle\neq 0$.
Then $\vhi\in\cE(L, t)$ with $t\in L^\ast$ a $G^\ast$-conjugate of $s$, by
\cite[Prop.~3.3.20]{GM20}. This in turn implies that $t$ is $N^\ast$-conjugate
to~$s$, since $N^\ast$ controls fusion in $L^\ast$ by \cite[Prop.~5.11]{Ma07}.
That is, such a character $\vhi$ lies in $\cE(L, s^{g^\ast})$ for some
$g^\ast\in N^\ast/L^\ast\cent_{N^\ast}(s)$.  It follows that
$\vhi(x)=\tht^g(x)\vhi(1)/\tht(1)$. Now, note that if
$g^\ast, h^\ast\in N^\ast/L^\ast\cent_{N^\ast}(s)_\la$ define the same coset
in $N^\ast/L^\ast\cent_{N^\ast}(s)$, then $\cE(L, s)^g=\cE(L, s)^h$.
Hence there is some $r'\in\QQ$ such that 
\[\sum_{\vhi\in\Irr(L)} \langle \chi, \RLG(\vhi)\rangle \vhi(x)
  =r'\sum_{g\in N/N_\tht}\sum_{\vhi\in\cE(L,s)^g}
     \langle\chi,\RLG(\vhi)\rangle\vhi(x)
  =\frac{r'}{\tht(1)}\sum_{g\in N/N_\tht} \tht^g(x)\ r_g\]
where $r_g:=\sum_{\vhi\in\cE(L, s)^g}\langle\chi,\RLG(\vhi)\rangle\vhi(1)$
for $g\in N/N_\tht$. It now suffices to see that $r_g=r_1$ for each
$g\in N/N_\tht$. Indeed, this follows since for $\vhi_1\in\cE(L,s)$, we have
$\vhi_1^g \in\cE(L,s)^g$ and $\langle\chi,\RLG(\vhi_1^g)\rangle\vhi_1^g(1)
 =\langle \chi, \RLG(\vhi_1)\rangle\vhi_1(1)$. 
\end{proof}

\begin{rem}   \label{rem:gendHC}
Note that since $\bL$ is a minimal $d$-split Levi subgroup of $\bG$, each
$\vhi\in\Irr(L)$ contributing to the sum in \eqref{eq:valchi} is
$d$-cuspidal. A conjecture of Cabanes--Enguehard (see the discussion after
\cite[Not.~1.11]{CE99}) posits that the $d$-cuspidal pairs under a given
character of $G$ must be $G$-conjugate. If this conjecture holds, we would
have that the $\vhi\in\Irr(L)$ contributing non-trivially to the sum
in~\eqref{eq:valchi} are always of the form $\vhi=\tht^g$ for some $g\in G$. In
the literature, it is sometimes said that \emph{generalised $d$-Harish-Chandra
theory} holds for $\cE(G,s)$ if moreover there is a bijection with the
characters of relative Weyl groups (see \cite[Def.~2.9]{KM13}).
We note that by \cite[Prop.~2.2.2]{Eng13}
and \cite[Thm~4.1]{CE99}, the conjecture of Cabanes--Enguehard holds (given
the properties of $\ell$ discussed in Section \ref{subsec:setting}) whenever
$s$ is an $\ell'$-element. Cabanes--Enguehard's conjecture has also been shown
to hold when $s$ is $\ell$-singular, under some assumptions on the prime~$\ell$
by Rossi \cite[Prop.~4.5, Cor.~4.11]{rossi}.
Further, when $d$ is a regular number and thus $\bL$ is an $F$-stable maximal
torus, we have $\langle \chi, \RLG(\vhi)\rangle\ne0$ if and only if
$(\bL, \vhi)$ is $G$-conjugate to $(\bL, \tht)$, if and only if
$(\bL, \vhi)=(\bL, \tht^g)$ for some $g\in N=\norm_G(\bL)$, using
\cite[Prop.~2.5.18 and Defs~2.5.14,~2.5.17]{GM20}. Further, in this case,
$\RLG(\vhi)=\RLG(\tht)$, by \cite[Cor.~2.2.10]{GM20}.

In particular, if we know that Cabanes--Enguehard's conjecture holds for
$\cE(G,s)$ (or, the more modest assertion that every $\vhi\in\Irr(L)$
contributing non-trivially to the sum in \eqref{eq:valchi} is of the form
$\vhi=\tht^g$ for $g\in G$), then we obtain an alternative proof for
Proposition~\ref{prop:nonreg} since then
\[\chi(x)=\tw{\ast}\RLG(\chi)(x)=r\sum_{g\in N/N_\tht} \tht^g(x)\]
for some integer $r\neq 0$, using \eqref{eq:valchi} and
Proposition~\ref{prop:dHCseries}. Hence $\QQ(\chi(x))=\QQ(\Om(\chi)(x))$ when
combined with \eqref{eq:local2}.
\end{rem}

\subsection{On the Strong Picky Conjecture}
Keep the situation of Setting \ref{subsec:setting}, with
$\chi\in\Irr_{\ell'}(G)\cap\cE(G,s)$ for $s\in G^\ast$ semisimple and
$\Om(\chi)\in\Irr_{\ell'}(N)$ corresponding to $(\tht,\eta)$ as in
Sections~\ref{subsec:local} and \ref{subsec:global}.

In some situations, such as when $\cent_{G^\ast}(s)\leq L^\ast$ or when $L$ is
a torus and $N_\tht=L$ (i.e., $\tht$ is in \emph{general position}), we have
$\pm\RLG(\tht)\in\Irr(G)$ (see \cite[Thm~3.3.22, Cor.~2.2.9]{GM20}). In this
case, we obtain the strong form of the conjecture:

\begin{prop}   \label{prop:genpos}
 Keep the setting of Section \ref{subsec:setting}. If
 $\chi=\pm \RLG(\tht)\in\Irr_{\ell'}(G)$, then~\ref{values} in
 Definition~\ref{def:strong} holds for $\chi$ at the prime $\ell$.

 That is, the Strong Picky Conjecture holds for the subset $\subG$ of
 characters of the form $\pm \RLG(\tht)\in\Irr_{\ell'}(G)$ for a minimal
 $d$-split Levi subgroup~$\bL$ and the restriction of $\Om$ to $\subG$.
\end{prop}

\begin{proof}
Let $x\in G$ be a picky $\ell$-element. Replacing $x$ by a suitable conjugate,
by Lemma \ref{lem:strab} we may assume that $x\in L$. Now we have
$\RLG(\tht)\otimes \St_\bG=\pm\Ind_{L}^G(\tht\otimes \St_{\bL})$ by
\cite[Cor.~3.4.11]{GM20}. Since $x$ is an $\ell$-element and
$\ell\nmid |\bZ(\bG)/\bZ^\circ(\bG)|$, its centraliser $\cent_\bG(x)$ is
connected (see, e.g.\ \cite[Exer.~20.16]{MT11}). Recall from
Theorem~\ref{thm:gunter} that the picky element $x$ satisfies $\cent_{G}(x)=L$.
Then from \cite[Prop.~3.4.10]{GM20}, we have
\[\pm\St_\bG(x)=|\cent_{G}(x)|_p=|L|_p.\]
From this, we see
\[\pm|L|_p\cdot \RLG(\tht)(x)
  =\Ind_{L}^G(\tht\otimes \St_{\bL})(x)
  =\frac{1}{|L|}\sum_{g\in G; x^g\in L} \tht(x^g)\St_\bL(x^g)
  =\frac{1}{|L|}\sum_{g\in N} \tht(x^g)\St_\bL(x^g),\]
where the last equality is from Lemma \ref{lem:strab}(e).

Again using \cite[Prop.~3.4.10]{GM20}, Theorem~\ref{thm:gunter}, and that $x$
is picky, for each $g\in N$, we have
$$\St_\bL(x^g)=\eps_\bL\eps_{\cent^\circ_\bL(x^g)} |\cent^\circ_\bL(x^g)^F|_p
  =\eps_\bL\eps_{\bL} |\cent_L(x^g)|_p=|L|_p$$
with signs $\eps_\bG\in\{\pm1\}$ as in \cite[Def.~2.2.11]{GM20}.
Then
$\RLG(\tht)(x)=\pm|L|^{-1}\sum_{g\in N} \tht(x^g)
  =\pm|L|^{-1}\sum_{g\in G; x^g\in L} \tht(x^g).$
That is, we see 
\begin{equation}\label{eq:RLG}
\RLG(\tht)(x)=\pm\Ind_L^G(\tht)(x)=\pm\Ind_L^N(\tht)(x).
\end{equation}

Recall that Conjecture~\ref{conj:AN}\ref{deg} holds for $G$ by
Lemma~\ref{lem:abelian} combined with \cite[Thm~7.8]{Ma07} and
\cite[Thm~1.1]{spath10}. The claim then follows from
Equations~\eqref{eq:local} and~\eqref{eq:RLG}, once we know that
$\weyl_G(\bL)_\theta$ (hence $\eta$) is trivial.

Now, the latter holds when the Mackey formula holds, as in
\cite[Cor.~9.3.1]{DM20}. This leaves only the cases $G=\tw{2}\type{E}_6(2)$ and
$G=\type{E}_8(2)$ when $\bL$ is not a torus (hence $d$ is not a regular number)
to consider, by \cite[Thm~3.3.7]{GM20}. In these cases, the Sylow
$\ell$-subgroups of $G$ are cyclic of prime order, and hence by
Lemma~\ref{lem:cyclic}, the result still holds. 
\end{proof}

We next obtain the Strong Picky Conjecture under suitable assumptions.

\begin{prop}   \label{prop:strongpickyreg}
 Keep the notation from Section \ref{subsec:setting}, and assume that (1) of
 Proposition \ref{prop:dHCseries} holds and that Cabanes--Enguehard's
 conjecture holds under $\chi$. (See Remark \ref{rem:gendHC}.) Then 
 \[\chi(x)=\pm\Om(\chi)(x).\]

 In particular, if $d$ is a regular number, then the Strong Picky Conjecture
 holds for $(G,\ell,x)$ for all picky $\ell$-elements $x$.
\end{prop}

\begin{proof}
From \eqref{eq:local} and the line before, we have
$\Om(\chi)(x)=\eta(1)\sum_{g\in N/N_\theta}\theta^g(x).$ Since (1) of
Proposition \ref{prop:dHCseries} holds, there is an $\wt{L}$-conjugate $\tht'$
of $\tht$ such that $\langle\chi,\RLG(\tht')\rangle\neq 0$. Arguing as in the
beginning of the proof of Proposition \ref{prop:nonreg}, we may assume instead
that $\langle\chi,\RLG(\tht)\rangle\neq 0$. On the other hand,
from~\eqref{eq:valchi} and the proof of Proposition \ref{prop:nonreg}, we know
then that
\[\chi(x)=\langle\chi, \RLG(\tht)\rangle \sum_{g\in N/N_\theta} \theta^g(x),\]
since we have assumed that every $\vhi\in\Irr(L)$ contributing to the sum
in~\eqref{eq:valchi} is conjugate to $\theta$.

We make the identification
$\weyl_{\cent_{G^\ast}(s)}(\cent^\circ_{\bL^\ast}(s))_\la
  \cong \weyl_G(\bL)_\tht$ as in~\eqref{eq:relweyl}. By combining
\cite[Thm~4.6(b)]{Ma07}, Brou{\'e}--Malle--Michel's Comparison Theorem
\cite[4.6.21]{GM20}, and \eqref{eq:IndRLG}, we have
\[\pm\eta(1)=\pm \langle\eta, \Ind_1^{\weyl_G(\bL)_\tht}(1)\rangle
  =\langle J_s^\bG(\chi), \Ind_{\cent_{G^\ast}^\circ(s)}^{\cent_{G^\ast}(s)}(R_{\cent_{\bL^\ast}^\circ(s)}^{\cent_{\bG^\ast}^\circ(s)}(\la))\rangle
\]
\[=\langle J_s^\bG(\chi), R_{\cent_{\bL^\ast}(s)}^{\cent_{\bG^\ast}(s)}(\Ind_{\cent_{L^\ast}^\circ(s)}^{\cent_{L^\ast}(s)}(\la))\rangle.\]

Now, recall from the proof of Proposition \ref{prop:dHCseries} that if
$\langle J_s^\bG(\chi), R_{\cent_{\bL^\ast}(s)}^{\cent_{\bG^\ast}(s)}(\hat\la)\rangle\neq0$ for some
$\hat\la\neq J_s^{\bL}(\tht)$ lying above $\la$, then $\hat\la=J_s^\bL(\tht')$
for some $\wt{L}$-conjugate $\tht'$ of $\tht$. But since (1) of
Proposition~\ref{prop:dHCseries} holds, we then have
$\langle \chi, \RLG(\tht')\rangle\neq0$. From our assumption that
Cabanes--Enguehard's conjecture holds under $\chi$, this means $\tht'$ is
further $N$-conjugate to $\tht$. That is, $\tht'=\tht^g=\tht^n$ for
$g\in\wt{L}$ and $n\in \norm_{G}(\bL)$. Then note that $s^{n^\ast}$ is
$L$-conjugate to $s$ (where $n^\ast\in N^\ast$ is associated to $n$ as before),
as $\tht$ and $\tht^n$ lie under the same character $\wt{\tht}$ of~$\wt{L}$.
Hence $n^\ast\in \cent_{N^\ast}(s)L^\ast$. Applying the equivariance of
Jordan decomposition from \cite[Thm~3.1]{CS13} to $\wt{L}$, we see that since
$\wt\tht^n$ corresponds to $\la^{n^\ast}$ but also to $\la$ (since
$J_s^\bL(\tht^n)$ lies over $\la$), it follows that in fact
$n^\ast\in \cent_{N^\ast}(s)_\la L^\ast$. Hence by the
isomorphism~\eqref{eq:relweyl}, we have $\tht^n=\tht$. That is, $\tht'=\tht$
and $\hat\la=J_s^\bL(\tht)$. Together with the computation above, this yields
\[\pm\eta(1)=\langle J_s^\bG(\chi), R_{\cent_{\bL^\ast}(s)}^{\cent_{\bG^\ast}(s)}(J_s^\bL(\tht))\rangle
 =\langle\chi, \RLG(\tht)\rangle,\]
where we again have used that (1) of Proposition \ref{prop:dHCseries} holds.
 
This gives 
\[\chi(x)=\pm\eta(1)\sum_{g\in N/N_\theta} \theta^g(x)=\pm\Om(\chi)(x),\]
as stated.

For the last statement, recall that, as noted in Remark~\ref{rem:regdHC},
(1) of Proposition~\ref{prop:dHCseries} holds when $d$ is a regular number
(i.e., when $\bL$ is a torus). Further, in that case generalised
$d$-Harish-Chandra theory holds by properties of Deligne--Lusztig characters
(see \cite[Cor.~2.2.10]{GM20}).
\end{proof}

\subsection{Suzuki and Ree groups}   \label{subsec:suzree}
We now assume that $F$ is a Steinberg, but not Frobenius, morphism. That is,
$G$ is a Suzuki or Ree group.

\begin{prop}\label{prop:suzree}
 The Picky+ Conjecture holds for any quasi-simple group $H$ such that
 $S=H/\bZ(H)$ is a Suzuki or Ree group ($\tw{2}\type{B}_2(2^{2f+1})$ with
 $f\geq 1$, $\tw{2}\type{G}_2(3^{2f+1})'$ with $f\geq 0$ or
 $\tw{2}\type{F}_4(2^{2f+1})'$ with $f\geq 0$) for all primes $\ell$ distinct
 from the defining prime. (That is, $\ell\neq 2$, $3$, resp. $2$.)
\end{prop}

\begin{proof}
For any $S$ as in the statement, the Schur multiplier of $S$ is trivial unless
$S=\tw{2}\type{B}_2(8)$, in which case it is a Klein four group. In the case
$S=\tw{2}\type{B}_2(8)$, the Sylow $\ell$-subgroups of~$S$, and hence of $H$,
are cyclic so this case is covered by Lemma~\ref{lem:cyclic}.
So, we may assume that $S$ is its own universal covering group. 

Further, we note that the case $S=\tw{2}\type{G}_2(3)'\cong \PSL_2(8)$ has been
handled above in Proposition \ref{prop:nonreg}, so we may assume $S$ is not
this group. If $S=\tw{2}\type{F}_4(2)'$, then $S\lhd G:=\tw{2}\type{F}_4(2)$
with index~2. In particular, $S$ and $G$ have the same picky $\ell$-elements.
(See \cite[Lemma~2.3]{Ma25}.) From here, the statement can be completed by
checking with \GAP. Hence, we may assume further that
$S\neq \tw{2}\type{F}_4(2)'$.

Then from now on, we may assume $S$ is of the form $G=\bG^F$ with $\bG$
simple of simply connected type and $F$ a Steinberg morphism. 
By \cite[Thm~5.11]{Ma25}, here an $\ell$-element $x\in G$ is picky if and only
if $\ell\geq 5$ and $x$ is a regular $\ell$-element. For $\tw{2}\type{B}_2(q^2)$
and $\tw{2}\type{G}_2(q^2)$, the Sylow $\ell$-subgroups are then cyclic, so
we are again done by Lemma~\ref{lem:cyclic}. 

Hence, we are left to consider the group $G=\tw{2}\type{F}_4(q^2)$, where
$q^2=2^{2f+1}$. Here the Sylow $\ell$-subgroups are again abelian, so
Lemma~\ref{lem:abelian} and Remark \ref{rem:picky+} apply. As noted in the
proof of \cite[Thm~5.1]{Malle14}, the conclusion of Lemma~\ref{lem:strab} still
holds here, taking $\bL$ to be a minimal $\Phi$-split Levi subgroup as defined
in loc.\ cit. The structure of $L=\bL^F$ (which is abelian) and
$N=\norm_G(P)=\norm_{G}(\bL)$ is given in \cite[Tab.~2]{Malle14}, and a
bijection analogous to $\Om$ above is given in \cite[Thm~8.5]{Ma07}. Note that
again $\cent_\bG(x)=\bL$. From here, the same considerations as in
Sections~\ref{subsec:setting}--\ref{subsec:global} yield the statement for
$G=\tw{2}\type{F}_4(q^2)$, where $q^2=2^{2f+1}$ with $f\geq 1$. (Indeed, the
extension property for $\Irr(L)$ with respect to $L\lhd N$ in
Section~\ref{subsec:local} is now given by \cite[Thm~D]{spath09}, and using
the map $\Om$ from \cite[Thm~8.5]{Ma07} we again obtain~\eqref{eq:local}
and~\eqref{eq:local2}. Further, \eqref{eq:valchi} still applies since
$\cent_\bG(x)=\bL$, still using \cite[Cor.~3.3.13]{GM20}. As $\bL$ is abelian
and $\bG=\wt\bG$ has connected centre, the analyses corresponding to those in
Proposition~\ref{prop:dHCseries} and~\ref{prop:nonreg} are simplified as in
Remarks~\ref{rem:regdHC} and~\ref{rem:gendHC}.)
\end{proof}

\section{Picky $2$-Elements}   \label{sec:picky 2}

In this section, we complete the determination of quasi-simple groups of Lie
type with picky semisimple $\ell$-elements started in \cite{Ma25} by
considering the remaining open case $\ell=2$. Throughout, $\bG$ is a simple
linear algebraic group of simply
connected type over an algebraically closed field of odd characteristic $p>2$,
$F:\bG\to\bG$ is a Frobenius map with respect to an $\FF_q$-structure and
$G:=\bG^F$. In order to determine picky 2-elements, we make use of the fact
that Sylow 2-subgroups of~$G$ are self-normalising most of the time (see
\cite[Cor.]{Ko05}):

\begin{prop}   \label{prop:syl2}
 Let $G$ be as above. Then a Sylow $2$-subgroup $P$ of $G$ is self-normalising
 except possibly when one of:
 \begin{enumerate}[\rm(1)]
  \item $G=\SL_n(\eps q)$ with $\eps\in\{\pm1\}$, $n\ge3$, where
   $\norm_G(P)/P\cong C_{(q-\eps)_{2'}}^{t-1}$;
  \item $G=\Sp_{2n}(q)$ with $n\ge1$ and $q\equiv3,5\pmod8$, where
   $\norm_G(P)/P\cong C_3^t$; or
  \item $G=\type{E}_6(\eps q)$ with $\eps\in\{\pm1\}$, where
   $\norm_G(P)/P\cong C_{(q-\eps)_{2'}}$.
 \end{enumerate}
 Here, $t$ in (1) and (2) denotes the number of non-zero digits in the $2$-adic
 expansion of $n$.
\end{prop}

\begin{lem}   \label{lem:torus}
 Let $G$ be as above and $x\in G$ be a picky $2$-element. Then $x$ lies in
 an $F$-stable maximal torus $\bT$ of $\bG$ with $|\bT^F|=(q-1)^a(q+1)^b$ for
 some $a,b\ge0$. Moreover, $x$ lies in a unique $F$-stable maximal torus
 of~$\bG$ except possibly when $G=\Sp_{2n}(3)$ for some $n\ge1$.
\end{lem}

\begin{proof}
Since $x$ is semisimple, it is contained in some $F$-stable maximal torus $\bT$
of~$\bG$. The order polynomial of $\bT$ is a product of cyclotomic polynomials.
Assume it is divisible by~$\Phi_d$ for some $d>2$. Since $q$ is odd, by
Zsigmondy's theorem there
exists a primitive prime divisor $r>d>2$ of $\Phi_d(q)$ and hence of $|\bT^F|$,
so of $|\cent_G(x)|$. But by Proposition~\ref{prop:syl2} the order of a Sylow
2-normaliser is only divisible by prime divisors of $6(q^2-1)$, so we obtain a
contradiction to \cite[Lemma~2.1]{Ma25}.

Now assume $x$ lies in at least two $F$-stable maximal tori of $\bG$. Then
$\cent_\bG(x)$ is not a torus, so has positive semisimple rank. In particular
$\cent_\bG(x)$ contains unipotent elements, and then so does $\cent_G(x)$. But
again by Proposition~\ref{prop:syl2} the order of a Sylow 2-normaliser is prime
to $p$ unless
$\bG=\Sp_{2n}$ and $p=3$. So assume $G=\Sp_{2n}(q)$ with $q=3^f$. As
$\cent_\bG(x)$ has semisimple rank at least~1, then $\cent_G(x)$ contains an
$\type{A}_1$-type subgroup, of order divisible by $(q^2-1)/2$. Now note that
$(q^2-1)/2$ is prime to~3, whence again using Proposition~\ref{prop:syl2},
\cite[Lemma~2.1]{Ma25} forces $q^2-1$ to be a 2-power, which only holds when
$q=3$.
\end{proof}

\begin{prop}   \label{prop:large rank}
 Assume that $G$ is one of $\SL_n(\eps q)$ with $n\ge6$, $\Sp_{2n}(q)$ with
 $n\ge6$, $\type{E}_6(\eps q)$, $\type{E}_7(q)$, or $\type{E}_8(q)$.
 Then $G$ possesses no picky $2$-element.
\end{prop}

\begin{proof}
For each group in the hypothesis, we will exhibit a subgroup $H\le G$ such that
\begin{enumerate}[(1)]
\item[(1)] $H$ contains representatives for all classes of maximal tori of $G$
 of order $(q-1)^a(q+1)^b$;
\item[(2)] $\norm_G(H)$ contains a Sylow $2$-subgroup of $G$; and
\item[(3)] the Sylow 2-subgroups of $\norm_G(H)/H$ are not normal.
\end{enumerate}
If such $H$ exists, then by Lemma~\ref{lem:torus} and (1) any class of picky
2-elements has a representative in $H$; but by (2) and (3) any such lies in
at least two Sylow 2-subgroups of $\norm_G(H)$, hence of~$G$, whence no picky
2-element can exist.
\par
To ensure (1), we let $H=\bH^F$ for a suitably chosen $F$-stable connected
reductive subgroup $\bH\le\bG$ containing a maximal torus of $\bG$ such that
the Weyl
group of $\bH$ contains representatives of all $F$-classes of elements of order
at most~2 in the Weyl group $W$ of $\bG$. Then, by \cite[Prop.~25.3]{MT11},
$H$ will contain the $F$-fixed points of representatives for all classes of
$F$-stable maximal tori of $\bG$ of orders $(q-1)^a(q+1)^b$.
Condition~(2) will be checked from the order formulas, and (3) by inspection.
\par
We consider the various types in turn. First assume $G=\SL_n(\eps q)$ with
$n\ge6$, and let $m:=\lfloor n/2\rfloor$. Then $W\cong\fS_n$ contains $m+1$
conjugacy classes of elements of order at most~2, all of which possess
representatives in the natural Young subgroup $\fS_2^m$, the Weyl group of the
standard Levi subgroup $H$ of type $\GL_2(\eps q)^m$ of $G$. Thus $H$
satisfies~(1).
Now $N:=\norm_G(H)$ contains a Sylow 2-subgroup of $G$, and $N/H\cong\fS_m$,
see \cite[Thm~4.10.6 and Tab.~4.10.6]{GLS}. Since the Sylow 2-subgroups of
$\fS_m$ are not normal for $m\ge3$, we are done.

For $G=\Sp_{2n}(q)$ write $n=2m+c$ with $c\in\{0,1\}$ and let $H$ be a natural
subgroup $\Sp_4(q)^m\times\Sp_2(q)^c$. Observe that the Weyl group of $H$
contains representatives of all involution classes in $W$; furthermore
$\norm_G(H)\ge\Sp_4(q)\wr\fS_m\times\Sp_2(q)^c$ satisfies~(2).
For $m\ge3$, so $n\ge6$, $\fS_m$ has non-normal Sylow 2-subgroups.

For $G=\type{E}_6(\eps q)$ let $H=(q-\eps)^2.\type{A}_1(q)^4$, the centraliser
of a maximal commuting product of fundamental $\SL_2$-subgroups (see
\cite[Tab.~4.10.6]{GLS}). Then $H$ contains representatives of all five
classes of maximal tori of $G$ of order $(q-\eps)^a(q+\eps)^b$, $2\le a\le6$,
and (2) is satisfied as $\norm_G(H)/H\cong\fS_4$. This also shows~(3).

For $G=\type{E}_7(q)$ we choose a subgroup $H=\type{A}_1(q)^3.\type{D}_4(q)$,
with $\norm_G(H)/H\cong\fS_3$ and for $G=\type{E}_8(q)$ a subgroup
$H=\type{D}_4(q)^2$, with $\norm_G(H)/H\cong\fS_3\times2$ (see
\cite[Tab.\ p.~168]{MT11}). By direct computation in \Chevie, these
satisfy~(1), and then~(2) and~(3) are also seen to hold.
\end{proof}

\begin{prop}   \label{prop:reduction}
 Assume $G$ contains a picky $2$-element. Then one of the following holds:
 \begin{enumerate}[\rm(1)]
  \item $q=3$;
  \item $G=\SL_n(\eps q)$ for $n\in\{3,5\}$ with $q\equiv-\eps\pmod4$;
  \item $G=\Sp_{2n}(5)$ with $n\le5$; or
  \item $G=\SL_2(q)$ or $\Sp_4(q)$.
 \end{enumerate}
\end{prop}

\begin{proof}
Let $x\in G$ be a picky 2-element. By Proposition~\ref{prop:large rank} we may
assume $G$ is not $\SL_n(\eps q)$ with $n\ge6$ or $\type{E}_6(\eps q)$.
By Lemma~\ref{lem:torus}, $x$ lies in some maximal torus $T$ of $G$ of order
$(q-1)^a(q+1)^b$. We first assume $q\equiv1\pmod4$. Then $m:=(q+1)/2\ge3$ is
odd and prime to $q-1$. By Proposition~\ref{prop:syl2} the order of the
normaliser of a Sylow 2-subgroup of $G$ is not divisible by~$m$ unless
$G=\SU_n(q)$ with $n$ not a 2-power, or $G=\Sp_{2n}(5)$.  So except for those
latter cases we conclude $|T|=(q-1)^n$ using \cite[Lemma~2.1]{Ma25}.
If on the other hand $q\equiv-1\pmod4$ and $q>3$ then $m:=(q-1)/2\ge3$ is odd
and prime to $q+1$, and so unless $G=\SL_n(q)$ with $n$ not a 2-power we
infer that $|T|=(q+1)^n$.
\par
Thus, unless we are in one of the exceptions in (1)--(3), for $\bT\le\bG$ an
$F$-stable maximal torus with $\bT^F=T$, by \cite[Prop.~5.20]{Ma07} the
normaliser $N:=\norm_G(\bT)$ contains a Sylow 2-subgroup of~$G$. We have thus
shown that $x\in T\unlhd N$ with $N/T=W$ the relative Weyl group of a Sylow
$d$-torus of $\bG$, where $d=1$ if $q\equiv1\pmod4$ and $d=2$ otherwise.
If $W$ has non-normal Sylow 2-subgroups then $x$ cannot be picky in $N$ and
hence neither in $G$. If $W$ has normal Sylow 2-subgroups, it is of type
$\type{A}_1$ or $\type{B}_2$ by \cite[Prop.~5.5]{Ma25}, so $G$ is
$\SL_2(q)$ or $\Sp_4(q)$ as in~(4).
\end{proof}

We now discuss the cases left open in the previous result. 

\begin{rem}
 The picky $\ell$-elements of a given finite permutation group or matrix group
 over a finite field can be determined effectively using the criterion in
 \cite[Cor.~2.10]{Ma25}, which we have implemented in the \GAP\ system
 \cite{GAP}. This was used to deal with some of the small cases below.
\end{rem}

\begin{lem}   \label{lem:SL3}
 Let $G=\SL_3(\eps q)$ with $q\equiv-\eps\pmod4$. Then $G$ contains picky
 $2$-elements if and only if $q$ is a Fermat or Mersenne prime, or $q=9$.
\end{lem}

\begin{proof}
By Lemma~\ref{lem:torus} any picky 2-element $x\in G$ lies in a unique maximal
torus $T$ of~$G$, of order $(q-1)^a(q+1)^b$. If $|T|=(q-\eps)^2$ then
$|\cent_G(x)|$ is divisible by $(q-\eps)_{2'}^2$ which is only
possible for $q=3$ by Proposition~\ref{prop:syl2}. Otherwise, $|T|=q^2-1$,
which is divisible by an odd prime not dividing $q-\eps$ unless $q+\eps$ is a
2-power, so $q$ is a Mersenne or Fermat prime, or $q=9$.
In the latter case, the normaliser of $T$ contains a Sylow 2-subgroup~$P$
of~$G$ and $\norm_G(P)=\norm_G(T)=T.2$ by Proposition~\ref{prop:syl2}. Now let
$x\in T$ be a generator of the cyclic Sylow 2-subgroup $T_2$ of $T$. Then
$T_2$ has index~2 in $P$, so is normal in any Sylow 2-subgroup of $P$
containing it. In particular, any such Sylow 2-subgroup is contained in
$\norm_G(T_2)=\norm_G(T)=\norm_G(P)$ and so equals $P$. Hence $x$ is picky.
\end{proof}

The next assertion is also shown in \cite[Thm~6.1]{MMM}; we provide a short
proof for completeness.

\begin{lem}   \label{lem:SL2}
 Let $G=\SL_2(q)$ with $q\ge5$ odd. Then $G$ contains a picky $2$-element if
 and only if $q$ is a Fermat or Mersenne prime or $q=9$.
\end{lem}

\begin{proof}
Any 2-element of $G$ is contained in a maximal torus of $G$ of order $q-\eps$
for $\eps\in\{\pm1\}$ with $q\equiv\eps\pmod4$.
First assume $q\equiv\eps\pmod8$. Then the Sylow 2-normalisers in $G$ are
dihedral of order $2(q-\eps)_2$. If $q$ is neither a Fermat or a Mersenne
prime nor $q=9$, then $(q-\eps)_2<q-\eps$ and thus there are no picky
2-elements by \cite[Lemma~2.1]{Ma25}.
Now assume $q\equiv4+\eps\pmod8$. In this case, the Sylow 2-normaliser is
isomorphic to $\SL_2(3)$. Since the largest element order in that group is~6,
we conclude $q-\eps\le6$ if $G$ has a picky 2-element. This forces $q=5$,
a Fermat prime.

Conversely, arguing similarly to Lemma~\ref{lem:SL3}, it is easy to see that
for $q\ge5$ a Fermat or Mersenne prime or $q=9$, any 2-element of $G$ of order
at least~4 lies in a unique Sylow 2-subgroup.
\end{proof} 

\begin{lem}   \label{lem:Sp4}
 Let $G=\Sp_4(q)$ with $q$ odd. Then $G$ contains picky $2$-elements if and
 only if $q\ne5,7$ is a Fermat or Mersenne prime.
\end{lem}

\begin{proof}
Since $\Sp_4(3)\cong2.\SU_4(2)$ the claim in this case follows from
\cite[Prop.~3.2]{Ma25}, while $\Sp_4(5)$, $\Sp_4(7)$ and $\Sp_4(9)$ have no
picky 2-element by explicit computation. So now assume $q\ge11$ and
$q\equiv\eps\pmod4$ with $\eps\in\{\pm1\}$. The maximal tori of $G=\Sp_4(q)$
have orders $(q\pm\eps)^2$, $q^2-1$ and $q^2+1$, so using
Proposition~\ref{prop:syl2} and \cite[Lemma~2.1]{Ma25} we check that any
picky 2-element $x\in G$ must be contained in a torus $T$ of order $(q-\eps)^2$
and $q-\eps$ necessarily is a 2-power.

Thus $q-\eps=2^f$ with $f\ge4$. The normaliser $N:=\norm_G(T)=T.W(\type{B}_2)$
of $T$ is the normaliser of a Sylow 2-subgroup of $G$. Now by explicit
computation all elements of $N$ lying in one of the seven
non-trivial cosets of $T$ either have some fourth or eighth root of unity as
eigenvalue, or their set of eigenvalues is closed under negation. Thus, for
$x\in T$ an element whose set of eigenvalues consists of 16th roots of unity
but is not closed under negation, the only $G$-conjugates of $x$ in $N$ are its
eight $N$-conjugates in~$T$, and thus $x$ is picky by \cite[Lemma~2.9]{Ma25}.
\end{proof}

We thus reach the following classification of groups with a picky 2-element:

\begin{thm}   \label{thm:picky 2}
 Let $G$ be as above. Then $G$ contains a picky $2$-element if and only if $G$
 is one of
 \begin{enumerate}[\rm(1)]
  \item $\SL_2(q)$ where $q$ is a Fermat or Mersenne prime or $q=9$;
  \item $\SL_3(q)$ with $q\equiv3\pmod4$ where $q$ is a Mersenne prime;
  \item $\SU_3(q)$ with $q\equiv1\pmod4$ where $q$ is a Fermat prime or $q=9$;
  \item $\Sp_4(q)$ where $q\ne5,7$ is a Fermat or Mersenne prime;
  \item $\SL_4(3),\SU_3(3),\SU_4(3),\Spin_7(3),\Sp_6(3),\Spin_8^\pm(3),
   \type{G}_2(3)$.
 \end{enumerate}
\end{thm}

\begin{proof}
We go through the cases in Proposition~\ref{prop:reduction}. If $q=3$ first
assume $G=\Sp_{2n}(3)$. Let $x\in G$ be a picky 2-element. If $x$ has order at
least~16 then its centraliser contains a cyclic subgroup $\GL_1(81)$ and thus
has order divisible by~5, contradicting Proposition~\ref{prop:syl2}.
Furthermore, if the eigenspace of $x$ for some 8th root of unity has dimension
at least~2 the centraliser contains a subgroup $\SL_2(9)$, again of order
divisible by~5. Thus, if $2n\ge8$ there is a subgroup
$H:=\Sp_4(3)\times\Sp_{2n-4}(3)\le G$ such that $x\in H$ and the
component~$x_1$ of $x$ in the factor $H_1=\Sp_4(3)$ has order at most~4.
By direct computation, the subnormaliser in $H_1$ of any such $x_1$ has order
divisible by~9, hence the same holds for the subnormaliser of $x$ in $H$ and
thus in $G$, whence $x$ cannot be picky in $G$ by Proposition~\ref{prop:syl2}.
We conclude that $2n<8$.

The same reasoning applies for
$\Spin_{2n+1}(3)$ and $\Spin_{2n}^\pm(3)$, using that none of $\SO_5(3)$ and
$\GO_6^\pm(3)$ possesses picky 2- or 4-elements, to conclude that $2n+1\le7$,
respectively $2n\le8$. Using Proposition~\ref{prop:large rank} for the
remaining types we then arrive at the groups listed in~(5), $\SL_3(3)$ and
$\Sp_4(3)$ in~(2) and~(4), and the groups
$$\SL_5(3), \SU_5(3), \tw3\type{D}_4(3), \type{F}_4(3).$$
By calculations in \GAP\ \cite{GAP}, none of the latter different from
$\type{F}_4(3)$ possesses picky 2-elements. For $G=\type{F}_4(3)$ one observes
that none of its maximal tori of order $(q-1)^a(q+1)^b$ does contain regular
elements (this can also be checked from the data at \cite{Lue}). Thus the
centraliser of any 2-element $x\in G$ contains unipotent elements, while the
Sylow 2-subgroups are self-normalising by Proposition~\ref{prop:syl2}, and so
$x$ cannot be picky. 
\par
Next let $G=\SL_5(\eps q)$ with $q\equiv-\eps\pmod4$. By Lemma~\ref{lem:torus}
any picky 2-element $x\in G$ lies in a unique maximal torus $T$ of $G$, of
order $(q-1)^a(q+1)^b$. Now $|T|$ is divisible by $(q-\eps)^2$, so
$|\cent_G(x)|$ is divisible by $(q-\eps)_{2'}^2$ which is only possible for
$q=3$ by Proposition~\ref{prop:syl2}, a case we already discussed. The case of
$\SL_3(\eps q)$ was dealt with in Lemma~\ref{lem:SL3}, leading to~(2) and~(3).
\par
The argument for $\Sp_{2n}(5)$, $n\ge3$, is completely analogous to the one
for $\Sp_{2n}(3)$ above, showing that $n>3$ is not possible. Again by
explicit computation, $\Sp_6(5)$ has no picky 2-element.
The case of $\SL_2(q)$ is discussed in Lemma~\ref{lem:SL2}, giving~(1).
Finally, $G=\Sp_4(q)$ leads to~(4) by Lemma~\ref{lem:Sp4}.
\end{proof}

The question of whether there are infinitely many groups $G$ as in the previous
statement hence hinges on the number theoretic question of existence or
non-existence of infinitely many Fermat or Mersenne primes.

\section{The Picky Conjecture for $\ell\le3$ and for exceptional covering groups}  \label{sec:23exceptcover}
We now complete the proofs of Theorems~\ref{thm:pickyodd}
and~\ref{thm:pickyeven} and derive Corollary~\ref{cor:quasisimple}.

\subsection{The case $\ell=2$}
To prove Theorem \ref{thm:pickyeven}, we first explicitly classify all picky
semisimple $2$-elements:

\begin{lem}   \label{lem:picky2classif}
 Let $G$ be the full covering group of one of the groups in
 Theorem~\ref{thm:picky 2}. Then the picky $2$-elements $x\in G$ are as follows:
 \begin{enumerate}[\rm(1)]
  \item in $\SL_2(q)$ with $3<q=2^n\pm1$: any $2$-element of order at least~$4$;
  \item in $\SL_3(q)$ with $q=2^n-1$: any $2$-element of order at least~$4$;
  \item in $\SU_3(q)$ with $q=2^n+1$: any $2$-element of order at least~$4$;
  \item in $\Sp_4(q)$ with $5,7\ne q=2^n\pm1$: certain regular $2$-elements of
   order at least~$8$; or
  \item $G,o(x)$ and $\cent_G(x)$ are as in Table~\ref{tab:picky 2}.
 \end{enumerate}
\end{lem}

\begin{table}[htb]
\caption{Picky 2-elements}   \label{tab:picky 2}
$$\begin{array}{c|cclc}
 G& o(x)& |\cent_G(x)|& |\Irr_{2^i}^x(G)|\cr
\hline
   \SL_4(3)& 8& 16\ (2\times)& 12\\
   \SU_3(3)& 8& 8\ (2\times)& 8\\
  3_i.\SU_4(3)& 8& 96\ (4\times)& 24,6,24\\
  6.\Om_7(3)& 8& 96\ (2\times)& 48,0,0,12\\
            & 8& 48& 48\\
   \Sp_6(3)& 8& 192\ (2\times)& 32, 24\\
           & 8& 32& 32\\
 \Om_8^+(3)& 8& 64\ (2\times)& 16,0,4,8\\
           & 8& 32& 16,0,4\\
 2.\Om_8^-(3)& 8& 64\ (2\times)& 32,0,0,8\\
  3.\mathrm{G}_2(3)& 8& 24\ (2\times)& 24\\
\end{array}$$
In the column headed $|\Irr_{2^i}^x(G)|$ we give the number of characters\\ in
$\Irr^x(G)$ with 2-part of the degree equal to $2^i$, for $i=0,1,2,\ldots$
\end{table}

\begin{proof}
The cases in (1)--(4) were discussed in Lemmas~\ref{lem:SL3}--\ref{lem:Sp4};
the cases in Table~\ref{tab:picky 2} are obtained by explicit computation in
these matrix groups.
\end{proof}

\begin{lem}   \label{lem:leftovers2}
 The Strong Picky Conjecture holds for all picky $2$-elements of the groups
 listed in Table \ref{tab:picky 2}, with a bijection preserving central
 characters.
\end{lem}

\begin{proof}
This can be checked in \GAP. The relevant numbers of non-vanishing
characters and the 2-parts of their degrees are displayed in
Table~\ref{tab:picky 2}.
\end{proof}

\begin{lem}   \label{lem:SL2conj}
 The Strong Picky Conjecture holds, with a bijection preserving central
 characters, for $G=\SL_2(q)$, $q$ odd, for any picky $2$-element $x\in G$.
\end{lem}

\begin{proof}
From Lemma \ref{lem:picky2classif} and the proof of Lemma \ref{lem:SL2},
if $G$ has a picky $2$-element $x\in P$ with $P\in\Syl_2(G)$, then either
$G=\SL_2(5)$ and $\norm_G(P)\cong \SL_2(3)$ or $q=2^n+\epsilon$ with $n\geq 3$
and $\eps\in\{\pm1\}$, and $\norm_G(P)\cong \Dih_{2^{n+1}}$. In either case,
recall that any $2$-element with $|x|\geq 4$ is picky. 

If $q=5$, we see directly from \GAP~that there is a unique class of picky
elements (of order 4) and both $G$ and $\norm_G(P)$ have four characters that
do not vanish on this class, all having odd degree and taking values in
$\{\pm1\}$ on $x$ for both groups and trivial on the centre.

Now let $q=2^n+\eps$ with $n\geq 3$, and let $x\in P$ be a $2$-element of
order larger than 2. From the well-known character tables of $G$ and
$\Dih_{2^{n+1}}$, we see that both groups have four characters
(all of odd degree) that take values $\pm1$ on $x$ and a family $\chi_k$ of
characters with $1\leq k< (q-\eps)/2$ that take values
$\pm(\zeta^k+\zeta^{-k})$ with $\zeta$ a root of unity satisfying
$|x|=|\zeta|$. The latter family has degree $2$ in $\norm_G(P)$ and degree
$q+\eps$ in $G$, so that $\chi_k(1)_2=2$ in both groups. Hence we see that
\ref{deg} and \ref{values} hold. Further, for both groups the odd-degree
characters are trivial on the centre and the families $\chi_k$ are trivial on
the centre exactly when $k$ is even. Hence the bijection further preserves
central characters.
\end{proof}

\begin{lem}   \label{lem:SL3conj}
 The Strong Picky Conjecture holds, with a bijection preserving central
 characters, for any picky $2$-element $x\in G=\SL_3(\eps q)$ with
 $\eps\in\{\pm1\}$ and $q$ odd.
\end{lem}

\begin{proof}
From Lemma \ref{lem:picky2classif}, we see that if picky $2$-elements exist,
then $q=2^n-\epsilon$ is a prime or $G=\SU_3(9)$, and then any $2$-element of
order larger than $2$ is picky. 

From the considerations in the proof of Lemma \ref{lem:SL3}, we have
$\norm_G(P)=T.2\cong C_{q^2-1}.2$. Here the order-2 action of $\norm_G(P)/T$
on the cyclic group $T$ is given by $a\mapsto a^{\epsilon q}$. Let $x$ have
order $|\zeta_2^a|$, where $\zeta_2$ is a primitive $(q^2-1)$th root of unity.
Now, from the character table in \Chevie\ \cite{Chev}, we see the characters
that are non-vanishing on $x$ are $1_G, \St_G$, which have odd degrees, take
values $\pm1$ on $x$, and are trivial on the centre; two families of characters
of odd degree indexed by $1\leq k<\frac{q-\eps}{2}$, which take values
$\pm\zeta_2^{ak(q+\eps)}$ on $x$ and lying above $\mu^k$ where $\mu$ is a
character generating the group $\Irr(\bZ(G))$; and one family of characters
indexed by $1\leq k<q^2-1$ with $(q+\eps)\nmid k$ whose degree is divisible
by~$2$ exactly once, which take values $\pm(\zeta_2^{ak}+\zeta_2^{-ak})$ on~$x$
and again lie over $\mu^k$. Then we see using Clifford theory that we have a
map satisfying \ref{deg} and \ref{values} and preserving central characters
by mapping the two extensions of $1_{C_{q^2-1}}$ to $\{1_G, \St_G\}$; the
extensions of the non-trivial characters of $T$ that are stable under
$\norm_G(P)/T$ to the two other families of odd-degree characters in
$\Irr^x(G)$; and the characters of $\norm_G(P)$ with reducible restriction
to~$T$ to the remaining family in $\Irr^x(G)$.
\end{proof}

\begin{lem}\label{lem:Sp4conj}
 The Strong Picky Conjecture holds, with a bijection preserving central
 characters, for $G=\Sp_4(q)$, $q$ odd, for any picky $2$-element $x\in G$.
\end{lem}

\begin{proof}
From Lemma \ref{lem:picky2classif} and the considerations in the proof of
Lemma~\ref{lem:Sp4}, we have $q=2^n+\eps$ with $\eps\in\{\pm1\}$ and $n\geq 4$,
and $x$ is a regular $2$-element lying in $T\cong C_{q-\eps}^2$. We further
have $P=\norm_G(P)=T.W(\type{B}_2)\cong (C_{q-\eps}.2)\wr C_2$. Here the
action of $C_{q-\eps}.2\setminus C_{q-\eps}$ is by inversion. From here, using
the character table for $G$ constructed by Srinivasan \cite{srinivasan} and
(somewhat tedious) Clifford theory arguments for $P$, we again obtain the
result.
\end{proof}

Lemmas \ref{lem:leftovers2}--\ref{lem:Sp4conj} complete the proof of
Theorem~\ref{thm:pickyeven}, when combined with Theorem~\ref{thm:picky 2} and
Lemma~\ref{lem:picky2classif}.

\subsection{The case $\ell=3$}
The simply connected quasi-simple groups of Lie type with non-abelian
Sylow 3-subgroups containing picky semisimple 3-elements were determined in
\cite[Thm~5.10]{Ma25}. By explicit computation in these groups in \GAP, one
obtains the following:

\begin{lem}   \label{lem:picky3classif}
 Let $G$ be the full covering group of one of the quasi-simple groups in
 \cite[Thm~5.10]{Ma25}. Then the picky $3$-elements $x\in G$, their orders and
 centralisers are as in Table~\ref{tab:picky 3}.
\end{lem}

\begin{table}[htb]
\caption{Picky 3-elements}   \label{tab:picky 3}
$\begin{array}{c|cclc}
 G& o(x)& |\cent_G(x)|& |\Irr_{3^i}^x(G)|\cr
\hline
  4_i.\SL_3(4)& 3& 36& 24\\
  \SU_3(8)& 9& 81\ (9\times)& 6,21\\
  2.\SU_4(2)& 9& 18\ (2\times)& 18\\
  \SU_5(2)& 9& 54\ (2\times), 27\ (2\times)& 27\\
  6.\SU_6(2)& 9& 324\ (6\times)& 36,72\\
            & 9& 54& 36\\
  \SU_7(2)& 9& 81& 54\\
  \SU_8(2)& 9& 486,243& 162\\
  2.\Sp_6(2)& 9& 18& 18\\
  \Sp_8(2)& 9& 54, 27& 27\\
  \Sp_{10}(2)& 9& 162\ (2\times)& 81\\
  2.\Om_8^+(2)& 9& 54\ (3\times)& 36\\
  \Om_{10}^+(2)& 9& 27& 27\\
  \Om_8^-(2)& 9& 9& 9\\
  \Om_{10}^-(2)& 9& 324, 162\ (2\times), 81& 54\\
  \Om_{12}^-(2)& 9& 162& 81\\
  \type{G}_2(8)& 9& 81\ (3\times)& 9,15\\
  \tw3\type{D}_4(2)& 9& 54\ (3\times)& 9,6\\
  2.\mathrm{F}_4(2)& 9& 108\ (2\times)& 54\\
\end{array}$
\end{table}

Again, in the column headed $|\Irr_{3^i}^x(G)|$ we give the number of
characters in $\Irr^x(G)$ with 3-part of the degree equal to $3^i$, for
$i=0,1,2,\ldots$.
Note that we do not consider the non-perfect group
$\type{G}_2(2)\cong\SU_3(3).2$.

\begin{lem}   \label{lem:leftovers3}
 The Strong Picky Conjecture holds, with a bijection preserving central
 characters, for all picky $3$-elements of the groups listed in
 Table~\ref{tab:picky 3}.
\end{lem}

\begin{proof}
All relevant character tables are contained in or can be computed by \GAP,
except for the tables of $\SU_8(2)$ and of $\type{G}_2(8)$. The generic
character table of $\type{G}_2(q)$ can be found in \cite{Chev}. From this,
by specialising $q$ to~$8$, one finds nine irreducible $3'$-characters and 15
characters $\chi$ with $\chi(1)_3=3$ for $G=\type{G}_2(8)$, as well as for
$\norm_G(P)$, satisfying~\ref{values} of Definition~\ref{def:strong}. Finally,
for $G=\SU_8(2)$ we use Lusztig's description of the irreducible characters.
The picky 3-elements lie in a maximal torus $T$ of order $\Phi_2^3\Phi_6^2$.
Any irreducible character of~$G$ taking non-zero value on one of those elements
has to lie in a Lusztig series $\cE(G,s)$ for which the centraliser of a
parametrising semisimple element $s\in G^*=\PGU_8(2)\cong G$ contains a torus
$T^*$ dual to $T$, so lies in $T^*$. Using the 8-dimensional matrix
representation of $G$ we determine all
such elements $s$ of order at most~3 up to $G$-conjugation. The values of
the characters in $\cE(G,s)$ on the two classes of picky 3-elements can then
be determined by the character formula. The values, up to sign, are collected
in Table~\ref{tab:SU8}. Here, the second column gives the number of elements
with the given centraliser, the next two columns show the multiplicities of
values (up to sign) on representatives $x_1,x_2$ of the two classes of picky
3-elements. Since the squares of these values add up to the centraliser order
of $x_i$ in $G$, all other irreducible characters of $G$ vanish on both $x_i$.
By direct computation, we obtain the same values up to sign (with
multiplicities) for the normaliser of a Sylow 3-subgroup of $G$.

The relevant numbers of non-vanishing characters for all groups in question are
displayed in Table~\ref{tab:picky 3}.
\end{proof}

\begin{table}[htb]
\caption{Non-vanishing characters of $\SU_8(2)$}   \label{tab:SU8}
$\begin{array}{c|ccc}
 \cent_G(s)& \#& |\chi(x_1)|& |\chi(x_2)|\\
\hline
                     \tw2\type{A}_7(q)& 1& 1\,(\ti6),2\,(\ti9),4\,(\ti3)& 1\,(\ti12), 2\,(\ti6)\\
              \Phi_2.\tw2\type{A}_6(q)& 2& 1\,(\ti6), 2\,(\ti3)& 1\,(\ti6), 2\,(\ti3)\\
       \Phi_2.\type{A}_1(q).\tw2\type{A}_5(q)& 2& 1\,(\ti6), 2\,(\ti9), 4\,(\ti3)& 1\,(\ti12), 2\,(\ti6)\\
            \Phi_2^2.\tw2\type{A}_5(q)& 1& 1\,(\ti6), 2\,(\ti3)& 1\,(\ti6), 2\,(\ti3)\\
   \Phi_2.\tw2\type{A}_2(q).\tw2\type{A}_4(q)& 2& 1\,(\ti9), 2\,(\ti9)& 1\,(\ti18)\\
            \Phi_2.\tw2\type{A}_3(q)^2& 1& 1\,(\ti9)& 1\,(\ti9)\\
 \Phi_2^2.\tw2\type{A}_2(q).\tw2\type{A}_3(q)& 2& 1\,(\ti9)& 1\,(\ti9)\\
   \Phi_2^2.\type{A}_1(q).\tw2\type{A}_2(q)^2& 1& 1\,(\ti9), 2\,(\ti9)& 1\,(\ti18)\\
\end{array}$
\end{table}

Lemma \ref{lem:leftovers3} completes the proof of Theorem~\ref{thm:pickyodd},
when combined with Propositions~\ref{prop:nonreg} and~\ref{prop:suzree}.

\subsection{Exceptional covering groups}
We now wish to prove the following:

\begin{prop}   \label{prop:exceptschur}
 Let $\ell$ be a prime and let $S$ be a simple group of Lie type defined
 in characteristic $p\neq\ell$, such that $S$ has an exceptional Schur
 multiplier. Then the Picky Conjecture holds for any quasi-simple group $H$
 such that $H/\bZ(H)=S$.
\end{prop}

\begin{proof}
By Lemma~\ref{lem:universalcoversuffices} we may assume $H$ is the universal
covering group of $S$, as long as the bijection for $H$ is chosen to preserve
central characters. There is some simple simply connected linear algebraic
group $\bG$ and Steinberg map $F$ such that $S=G/\bZ(G)$ for $G:=\bG^F$,
and we can write $G=H/Z$ for some $Z\leq \bZ(H)$. By \cite[Lemma~2.4]{Ma25},
the picky $\ell$-elements of $H$ are preimages of picky $\ell$-elements of~$G$
under the natural map.

First assume that $\ell\geq 5$. By \cite[Tab.~24.3]{MT11} we note that $|Z|$
is not divisible by $\ell$ so we know again using \cite[Thm~5.10]{Ma25} that
$P\in\Syl_\ell(H)$ is abelian if $H$ has picky $\ell$-elements. By
Lemma~\ref{lem:cyclic} we may also assume that Sylow $\ell$-subgroups of~$H$
and hence of $S$ are non-cyclic. This only leaves the possibilities that
$\ell=5$ and $S=\Om_8^+(2),\type{G}_2(4),\type{F}_4(2),\tw2\type{E}_6(2)$,
or $\ell=7$ and $S=\type{F}_4(2),\tw2\type{E}_6(2)$. Now, by the known
character tables, in any of these groups the centralisers of $\ell$-elements
are strictly bigger than the centraliser of an ambient Sylow $d$-torus, so by
\cite[Thm~5.3]{Ma25} none of them has picky $\ell$-elements, and so neither
has $H$.

For $\ell=2,3$ the validity of the Picky Conjecture was checked in
Lemmas~\ref{lem:leftovers2} and~\ref{lem:leftovers3} for all groups in
question apart from central quotients of $2.\SL_4(2)$ and $6.\tw2\type{E}_6(2)$,
for $\ell=3$. For $2.\SL_4(2)$ verifying the conjecture is an easy exercise
using \GAP, while $\tw2\type{E}_6(2)$ and hence $6.\tw2\type{E}_6(2)$ does not
possess picky 3-elements by inspection of the character table.
\end{proof}

Finally, we obtain Corollary \ref{cor:quasisimple}:

\begin{proof}[Proof of Corollary \ref{cor:quasisimple}]
By Lemma~\ref{lem:universalcoversuffices} we may assume $H$ is the full
covering group of the simple group of Lie type $S$ in question, as long as we
also know that the bijection used for $H$ can be chosen to preserve central
characters.
From Propositions \ref{prop:suzree} and \ref{prop:exceptschur}, we may assume
that $S$ is not a Suzuki or Ree group and that $S$ has a non-exceptional Schur
multiplier. In this case, the universal covering group is of the form
$G=\bG^F$ with $\bG$ simple of simply connected type and $F$ a Frobenius
morphism.  If $\ell=3$ and the Sylow $\ell$-subgroups of $G$ are non-abelian,
then by Lemma~\ref{lem:picky3classif} the relevant cases are listed in
Table~\ref{tab:picky 3}, and satisfy the Picky+ Conjecture, with preservation
of central characters, by Lemma~\ref{lem:leftovers3}. If $\ell=2$ the picky
elements for $G$ are collected in Lemma~\ref{lem:picky2classif} and they are
shown to satisfy the Picky Conjecture, with preservation of central characters,
in Lemmas~\ref{lem:leftovers2}--\ref{lem:Sp4conj}.

Otherwise, the Sylow $\ell$-subgroups of $G$ are abelian, and the statement
follows from the proof of Proposition \ref{prop:nonreg}, since $\Om$ preserves
central characters as noted after the statement of Theorem~\ref{thm:gunter}.
\end{proof}


\end{document}